\documentclass[11pt]{article}
\usepackage[T2A]{fontenc}
\usepackage[cp1251]{inputenc}

\usepackage[T2A]{fontenc}
\usepackage[cp1251]{inputenc}
\usepackage{amsfonts}
\usepackage{enumerate}
\usepackage{amssymb}
\usepackage{amsmath}

\begin{document}

\begin{center}
\LARGE\textbf{The Heyde Theorem on a-adic Solenoids}
\end{center}

\begin{center}
\textbf{M. Myronyuk}
\end{center}

\begin{center}
\textit{Mathematical Division, Institute for Low Temperature Physics
and Engineering,}

\textit{47, Lenin Ave, Kharkov, 61103, Ukraine}

E-mail address: myronyuk@ilt.kharkov.ua
\end{center}

\begin{abstract}
We prove the following analogue of the Heyde theorem for ${\bf
a}$-adic solenoids. Let $ \xi_1$, $\xi_2$ be independent random
variables taking values in an ${\bf a}$-adic solenoid $ \Sigma_{\bf
a}$ and with distributions $\mu_1$, $\mu_2$.  Let $\alpha_j,
\beta_j$ be topological automorphisms of $\Sigma_{\bf a}$ such that
$\beta_1\alpha^{-1}_1 \pm \beta_2\alpha^{-1}_2$ are topological
automorphisms of $\Sigma_{\bf a}$ too. Assuming that the conditional
distribution of the linear form $L_2=\beta_1\xi_1 + \beta_2\xi_2$
given $L_1=\alpha_1\xi_1 + \alpha_2\xi_2$ is symmetric, we describe
possible distributions $\mu_1$, $\mu_2$.
\end{abstract}



\section{Introduction}
\label{}

Many studies have been devoted to characterizing Gaussian
distributions on the real line. Specifically, in 1970 Heyde proved
the following theorem, which characterizes a Gaussian distribution
by the symmetry of the conditional distribution of one linear form
given another.

\medskip

\noindent\textbf{The Heyde theorem} [{Heyde (\cite{Heyde}; see also
\cite[Section 13.4.1]{KagLinRao})}]. \textit{Let $\xi_1,..., \xi_n$,
$n\ge 2,$ be independent random variables, $\alpha_j, \beta_j$ be
nonzero constants such that $\beta_i\alpha_i^{-1} \pm
\beta_j\alpha_j^{-1} \ne 0$ whenever $i \ne j$. If the conditional
distribution of the linear form $L_2 = \beta_1\xi_1 + \cdots +
\beta_n\xi_n$ given $L_1 = \alpha_1\xi_1 + \cdots + \alpha_n\xi_n$
is symmetric then all random variables $\xi_j$ are Gaussian.}

\medskip

In recent years, a great deal of attention has been focused upon
generalizing the classical characterization theorems to random
variables with values in locally compact Abelian groups (see e.g.
\cite{Fe3}-\cite{Fe-He-3}, \cite{Fe-solenoid-2009}-\cite{FeMy},
\cite{MiFe-He-1}, \cite{My-He-2009}; see also \cite{Fe-book2}, where
one can find additional references). The articles
\cite{Fe-He-1}--\cite{Fe-He-3}, \cite{MiFe-He-1}, \cite{My-He-2009}
(see also \cite[Chapter VI]{Fe-book2}) were devoted to finding
group-theoretic analogues of the Heyde theorem. This article
continues this research.

Let $X$ be a second countable locally compact Abelian group, ${\rm
Aut}(X)$ be the group of topological automorphisms of $X$. Let
$\xi_j, \ j = 1, 2, \dots, n, \ n \ge 2,$ be independent random
variables taking values in $X$ and with distributions $\mu_j$. Let
$\alpha_j, \beta_j \in {\rm Aut}(X)$ such that $\beta_i\alpha_i^{-1}
\pm \beta_j\alpha_j^{-1} \in {\rm Aut}(X)$ whenever $i\neq j$.
Define the linear forms $L_1=\alpha_1\xi_1 + \cdots + \alpha_n\xi_n$
and $L_2=\beta_1\xi_1 + \cdots + \beta_n\xi_n$.

We formulate the following problem.

\medskip

\noindent\textbf{Problem 1.} \textit{Assume that conditional
distribution of $L_2$ given $L_1$ is symmetric. Describe the
possible distributions $\mu_j$.}

\medskip

Problem 1 was solved for the class of finite Abelian groups in
\cite{Fe-He-1}, \cite{MiFe-He-1} and then for the class of countable
discrete Abelian groups in \cite{Fe-He-3}, \cite{My-He-2009}.
Problem 1 for \textbf{a}-adic solenoids was formulated in the book
\cite{Fe-book2}. In this article we solve this problem.

\textbf{a}-adic solenoids are important examples of connected
Abelian groups. We note that if $X$ is a connected Abelian group of
dimension one then $X$ is topologically isomorphic either the real
line $\mathbb{R}$, or the circle group $\mathbb{T}$, or an
\textbf{a}-adic solenoid $\Sigma_{\mathbf{a}}$. Problem 1 was solved
for the case $X=\mathbb{R}$ by Heyde. Problem 1 cannot be formulated
for the case $X=\mathbb{T}$ because there no exist topological
automorphisms $\alpha_j, \beta_j$ such that $\beta_i\alpha_i^{-1}
\pm \beta_j\alpha_j^{-1} \in {\rm Aut}(X)$ whenever $i\neq j$. In
this article we solve Problem 1 (Theorem 1) for \textbf{a}-adic
solenoids $\Sigma_{\mathbf{a}}$. It turns out that the answer
depends on topological automorphisms $\alpha_j, \beta_j$. Note that
it follows from \cite{Fe-He-2} (see also \cite[\S16.2]{Fe-book2})
that under the condition that the characteristic functions of the
distributions $\mu_j$ do not vanish the symmetry of the conditional
distribution of $L_2$ given $L_1$ implies that $\mu_j$ are Gaussian.

\section{Notation and definitions}

Let $X$ be a locally compact  Abelian group, $Y=X^\ast$ be its
character group, and  $(x,y)$ be the value of a character $y \in Y$
at an element $x \in X$. Let $K$ be a subgroup of $X$. Denote by
$A(Y, K) = \{y \in Y: (x, y) = 1$ for all $x \in K \}$ the
annihilator of $K$. If $\delta : X \mapsto X$ is a continuous
homomorphism, then the adjoint homomorphism $\widetilde \delta : Y
\mapsto Y$ is defined by the formula $(x, \widetilde \delta y) =
(\delta x, y)$ for all $x \in X, \ y \in Y$. We note that $\delta
\in {\rm Aut}(X)$ if and only if $\widetilde\delta \in {\rm
Aut}(Y)$. For each integer $n$, $n \ne 0,$ let $f_n : X \mapsto X$
be the homomorphism $f_n x=nx.$ Set $X^{(n)} = f_n(X)$, $X_{(n)}=Ker
f_n$. It is clear that the adjoint homomorphism
${\widetilde{f_n}:Y\mapsto Y}$ is the mapping $\widetilde{f_n}
y=ny$. Denote by ${\mathbb{R}}$ the additive group of real numbers,
by ${\mathbb{Z}}$  the group of integers, by $\mathbb{Q}$ the group
of rational numbers considering in the discrete topology, by
${\mathbb{Z}}(n)$  the finite cyclic group of order $n$. For a fixed
prime $p$ denote by ${\mathbb{Z}}(p^{\infty})$ the set of rational
numbers of the form $\{ k/p^n : k = 0, 1,..., p^n-1, n = 0, 1,...\}$
and define the operation in ${\mathbb{Z}}(p^{\infty})$ as addition
modulo 1. Then ${\mathbb{Z}}(p^{\infty})$ is transformed into an
Abelian group, which we consider in the discrete topology. Denote by
${\rm Aut}(X)$ the group of topological automorphisms of the group
$X$.

Put  ${\bf a}= (a_0, a_1,\dots)$, where all $a_j \in {\mathbb{Z}}, \
a_j > 1$. First we recall the definition of the group of ${\bf
a}$-adic integers $\Delta_{{\bf a}}$ \cite[(10.2)]{HeRo1}. As a set
$\Delta_{{\bf a}}$ coincides with the Cartesian product
$\mathop{\mbox{\rm\bf P}}\limits_{n=0}^\infty\{0,1,\dots ,a_n-1\}$.
For $\mathbf{x}=(x_0, x_1, x_2,\dots), \ \mathbf{y}=(y_0, y_1,
y_2,\dots)\in \Delta_{{\bf a}}$ let
$\mathbf{z}=\mathbf{x}+\mathbf{y}$ be define as follows. Let
$x_0+y_0=t_0a_0+z_0$, where $z_0\in\{0, 1,\dots ,a_0-1\}, \
t_0\in\{0, 1\}$. Assume that the numbers $z_0,z_1,\dots ,z_k;\ t_0,
t_1,\dots ,t_k$ have been already determined. Let us put then
$x_{k+1}+y_{k+1}+t_k=t_{k+1}a_{k+1}+z_{k+1}$, where
$z_{k+1}\in\{0,1,\dots ,a_{k+1}-1\},\ t_{k+1}\in\{0, 1\}$. This
defines by induction a sequence $\mathbf{z}=(z_0, z_1, z_2,\dots)$.
The set $\Delta_{{\bf a}}$ with the addition defined above is an
Abelian group, whose neutral element is the sequence in
$\Delta_{{\bf a}}$ that is identically zero. Consider $\Delta_{{\bf
a}}$ in the product topology. The obtained group is called the ${\bf
a}$-adic integers. If all of the integers $a_j$ are equal to some
fixed prime integer $p$, we write $\Delta_p$ instead of
$\Delta_{{\bf a}}$, and call this object the group of $p$-adic
integers. Note that $\Delta_p^*\approx {\mathbb{Z}}(p^{\infty})$
(see \cite[(25.2)]{HeRo1}).

Consider the group $\mathbb{R}\times\Delta_{{\bf a}}$. Let $B$ be
the subgroup of the group $\mathbb{R}\times\Delta_{{\bf a}}$ of the
form $B=\{(n,n\mathbf{u})\}_{n=-\infty}^{\infty}$, where
$\mathbf{u}=(1, 0,\dots,0,\dots)$. The factor-group $\Sigma_{{\bf
a}}=(\mathbb{R}\times\Delta_{{\bf a}})/B$ is called the ${\bf
a}$-adic solenoid. The group $\Sigma_{{\bf a}}$ is compact and
connected, and has dimension one ([\cite[(10.12), (10.13),
(24.28)]{HeRo1}]). The character group of the group $\Sigma_{{\bf
a}}$ is topologically isomorphic to the subgroup $H_{\bf a}\subset
\mathbb{Q}$ of the form
$$H_{\bf a}= \left\{{m \over a_0a_1 \dots a_n} :
\quad n = 0, 1,\dots; \ m \in {\mathbb{Z}} \right\}.
$$
We will assume without loss of generality that if $X=\Sigma_{{\bf
a}}$ then $Y=X^*=H_{\bf a}$.

Let $Y$ be an Abelian group, $f(y)$ be a function on $Y$, and $h\in
Y.$ Denote by $\Delta_h$ the finite difference operator
    $$\Delta_h f(y)=f(y+h)-f(y).$$
A function $f(y)$ on $Y$ is called a polynomial if
    $$\Delta_{h}^{n+1}f(y)=0$$
for some $n$ and for all $y,h \in Y$. If $Y$ is a subgroup of
$\mathbb{Q}$ then this definition of a polynomial coincides with the
classical one.

Let ${M^1}(X)$ be the convolution semigroup of probability
distributions on $X$,  $\hat \mu(y) = \int_X (x, y) d\mu(x)$ be the
characteristic function of a distribution $\mu \in {M^1}(X)$, and
$\sigma(\mu)$ be the support of $\mu$. If $K$ is a closed subgroup
of $X$ and $\sigma(\mu) \subset K$, then $\hat\mu(y+h) = \hat\mu(y)$
for all $y \in Y, \ h \in A(Y, K)$. If $E$ is a closed subgroup of
$Y$ and $\hat \mu(y)=1$ for $y \in E$, then $\hat\mu(y+h) =
\hat\mu(y)$ for all $y \in Y, \ h \in E$ and $\sigma(\mu) \subset
A(X, E)$. For $\mu \in {M^1}(X)$ we define the distribution $\bar
\mu \in M^1(X)$ by the rule $\bar \mu(B) = \mu(-B)$ for all Borel
sets $B \subset X$. Observe that $\hat {\bar \mu}(y) =
\overline{\hat \mu(y)}$.

A distribution $\gamma \in {M^1}(X)$ is called Gaussian (\cite[\S
4.6]{Pa}) if its characteristic function can be represented in the
form
$$ \hat\gamma(y)= (x,y)\exp\{-\varphi(y)\},$$
where $x \in X$ and $\varphi(y)$ is a continuous nonnegative
function satisfying the equation
\begin{equation}\label{fe1}
    \varphi(u+v)+\varphi(u-v)=2[\varphi(u)+
\varphi(v)], \quad u, \ v \in Y.
\end{equation}
 Denote by $\Gamma(X)$ the set
of Gaussian distributions on $X$. It is easy to see that any
nonnegative function $\varphi(y)$ on the group $H_{\bf a}$
satisfying equation (\ref{fe1}) is of the form $\varphi(y)=\lambda
y^2$, where $\lambda\geq 0$, $y\in H_{\bf a}$. It is well-known that
a support of a Gaussian distribution on a locally compact Abelian
group $X$ is a coset of a connected subgroup of $X$. Thus if
$\gamma$ is a non degenerate Gaussian distribution on
$X=\Sigma_{\mathbf{a}}$ then $\sigma(\gamma)=X$.

Denote by $I(X)$ the set of idempotent distributions on $X$, i.e.
the set of shifts of Haar distributions $m_K$ of compact subgroups
$K$ of the group $X$. Observe that the characteristic function of
the Haar distribution $m_K$ is of the form

\begin{equation}\label{def1}
    \widehat m_K(y) =
\left\{%
\begin{array}{ll}
    1, & \hbox{$y \in A(Y,K)$;} \\
    0, & \hbox{$y\not\in A(Y,K)$.} \\
\end{array}%
\right.
\end{equation}

We note that if a distribution $\mu \in \Gamma(X)*I(X)$, i.e.
$\mu=\gamma*m_K$, where $\gamma\in \Gamma(X)$, then $\mu$ is
invariant with respect to a compact subgroup $K \subset X$ and under
the natural homomorphism $X \mapsto X/K \ \mu$ induces a Gaussian
distribution on the factor group $X/K$.

\section{The Heyde theorem (the general case)}

Let $ \xi_1, \ \xi_2$ be independent random variables with values in
the group $X = \Sigma_{\bf a}$ and distributions $\mu_1, \ \mu_2$.
Consider the linear forms $L_1=\alpha_1\xi_1 + \alpha_2\xi_2$ and
$L_2=\beta_1\xi_1 + \beta_2\xi_2$, where $\alpha_j, \beta_j \in {\rm
Aut}(X)$ and $\beta_1\alpha^{-1}_1 \pm \beta_2\alpha^{-1}_2 \in {\rm
Aut}(X)$. Assume that the conditional distribution of linear form
$L_2$ given $L_1$ is symmetric. Taking into consideration new
independent random variables $\xi'_j = \alpha_j\xi_j$ we reduce the
study of the distributions $\mu_j$ on $X$ to the case when $L_1 =
\xi_1 + \xi_2$ and $L_2=\delta_1 \xi_1  + \delta_2 \xi_2 $, where
$\delta_j \in {\rm Aut}(X)$ and $\delta_1 \pm \delta_2 \in {\rm
Aut}(X)$. Note that any topological automorphism $\delta$ of the
group $X$ is of the form
$$ \delta
= f_p f_q^{-1}
$$
for some relatively prime $p$ and $q$, where $f_p, f_q \in {\rm
Aut}(X)$. Note that for any ${\delta \in {\rm Aut}(X)}$ the
conditional distribution of the linear form $L_2$ given $L_1$ is
symmetric if and only if the conditional distribution of the linear
form $\delta L_2$ given $L_1$ is symmetric. Hence without loss of
generality, we can assume from the beginning that $L_1 = \xi_1 +
\xi_2$, $L_2 = p \xi_1 +  q \xi_2$, where $p, q \in {\mathbb{Z}}$,
$p q \ne 0$, $p$ and $q$ are relatively prime, $f_p, f_q, f_{p\pm q}
\in {\rm Aut}(X)$. Now we formulate the main result of this article.

\medskip

\noindent\textbf{Theorem 1.} \textit{Let $X=\Sigma_{\mathbf{a}}$.
Assume that $f_p, f_q, f_{p\pm q} \in {\rm Aut}(X)$, $p$ and $q$ are
relatively prime. The following statements hold:}
\begin{enumerate}
    \item \textit{Assume that $pq=-3$. Let $\xi_1,\xi_2$ be
    independent random variables with values in
    $X$ and distributions $\mu_1,\mu_2$. If the
    conditional distribution of the linear form
    $L_2=p\xi_1+q\xi_2$ given $L_1=\xi_1+\xi_2$
    is symmetric then at least one distribution $\mu_j\in \Gamma(X)*I(X)$.}
    \item \textit{Assume that $pq\neq -3$. Then there exist
    independent random variables $\xi_1,\xi_2$ with values
    in $X$ and distributions $\mu_1,\mu_2$ such that the
    conditional distribution of the linear form $L_2=p\xi_1+q\xi_2$
    given $L_1=\xi_1+\xi_2$ is symmetric and the distributions
    $\mu_j\not\in \Gamma(X)*I(X)$, $j=1,2$.}
\end{enumerate}

\medskip

Theorem 1 can be regarded as a group analogue of the Heyde theorem
for \textbf{a}-adic solenoids. To prove Theorem 1 we need some
lemmas.

\medskip

\noindent\textbf{Lemma 1.} \textit{Let $X$ be a locally compact
second countable Abelian group. Let $\xi_1, \ \xi_2$ be independent
random variables with values in $X$ and distributions $\mu_1,
\mu_2$. Consider the linear forms $L_1=\alpha_1\xi_1 +
\alpha_2\xi_2$ and $L_2=\beta_1\xi_1 + \beta_2\xi_2$, where
$\alpha_j, \beta_j$ are continuous homomorphisms of the group $X$.
The conditional distribution of the linear form $L_2$ given $L_1$ is
symmetric if and only if the characteristic functions of the
distributions $\mu_j$ satisfy the equation}
\begin{equation}\label{l0}
    \hat\mu_1(\tilde{\alpha_1}u+\tilde{\beta_1}v)\hat\mu_2(\tilde{\alpha_2}u+\tilde{\beta_2}v)=
    \hat\mu_1(\tilde{\alpha_1}u-\tilde{\beta_1}v)\hat\mu_2(\tilde{\alpha_2}u-\tilde{\beta_2}v),\quad u,v \in Y.
\end{equation}

\medskip

Lemma 1 was proved in \cite[\S 16.1]{Fe-book2} in the case where
$\alpha_j, \beta_j \in {\rm Aut}(X)$. This proof is valid for
arbitrary continuous homomorphisms $\alpha_j, \beta_j$ of the group
$X$.

\medskip

\noindent\textbf{Lemma 2.} \textit{Let either $|q|=2$ or $q=4m+3$,
where $m$ is some integer. Let $X=\Delta_2$. Then there exist
independent identically distributed random variables $\xi_1,\xi_2$
with values in $X$ and distribution $\mu\not\in I(X)$ such that the
conditional distribution of the linear form $L_2=\xi_1+q\xi_2$ given
$L_1=\xi_1+\xi_2$ is symmetric.}

\medskip

\textbf{Proof.} Since $X=\Delta_2$, we have $Y\approx
\mathbb{Z}(2^{\infty})$. To avoid introducing new notation we will
assume that $Y=\mathbb{Z}(2^{\infty})$.

Let $g_0(y)$ be an arbitrary characteristic functions on $Y_{(2)}$.
Set

$$g(y)=
\left\{%
\begin{array}{ll}
    g_0(y), & \hbox{$y\in Y_{(2)}$;} \\
    0, & \hbox{$y\not\in Y_{(2)}$.} \\
\end{array}%
\right.$$

The function $g(y)$ is a positive definite function on $Y$ (\cite[\S
32]{HeRo2}]). By the Bochner theorem there exists a distribution
$\mu \in {M}^1(X)$ such that $\widehat\mu(y) = g(y)$. It is clear
that $g_0(y)$ can be chosen in such a way that $\mu\not\in I(X)$.

Let $\xi_1,\xi_2$ be independent identically distributed random
variables with values in $X$ and  distribution $\mu$. We check that
the conditional distribution of the linear form $L_2=\xi_1+q\xi_2$
given $L_1=\xi_1+\xi_2$ is symmetric. By Lemma~\ref{Lemma-1} it
suffices to show that the characteristic function $\hat\mu(y)$
satisfies equation (\ref{l0}) which takes the form

\begin{equation}\label{l1}
    \hat\mu(u+v)\hat\mu(u+qv)=\hat\mu(u-v)\hat\mu(u-qv),\quad u,v \in Y.
\end{equation}

Since $Y_{(2)}\approx \mathbb{Z}(2)$, it is clear that if $u, v \in
Y_{(2)}$ then equation (\ref{l1}) is an equlity.

If either $u\in Y_{(2)}, v\not\in Y_{(2)}$ or $u\not\in Y_{(2)},
v\in Y_{(2)}$ then $u\pm v\not\in Y_{(2)}$. Hence
$\hat\mu(u+v)=\hat\mu(u-v)=0$ and equation (\ref{l1}) is an
equality.

Let $u, v\not\in Y_{(2)}$. Suppose that the left-hand side of
equation (\ref{l1}) does not vanish. Then

\begin{equation}\label{l1.1}
    u+v\in Y_{(2)}, \quad u+qv\in Y_{(2)}.
\end{equation}

Let $q=2$. It follows from (\ref{l1.1}) that $v \in Y_{(2)}$,
contrary to the choice of $v$. Hence the left-hand side of equation
(\ref{l1}) is equal to zero. Similarly, we prove that the right-hand
side of equation (\ref{l1}) is equal to zero.

Let $q=-2$. It follows from (\ref{l1.1}) that $3v \in Y_{(2)}$.
Since $f_{3}\in\mathrm{Aut}(Y)$ and $Y_{(2)}$ is a characteristic
subgroup, we have $v\in Y_{(2)}$, contrary to the choice of $v$.
Hence the left-hand side of equation (\ref{l1}) is equal to zero.
Similarly, we prove that the right-hand side of equation (\ref{l1})
is equal to zero.

Let $q=4m+3$. It follows from (\ref{l1.1}) that $(q-1)v \in
Y_{(2)}$. Since $(q-1)= 2(2m+1)$ and $f_{2m+1}\in\mathrm{Aut}(Y)$,
we have $2v\in Y_{(2)}$. Hence $v$ is an element of order 4. So,
$qv=-v$. It follows from this that equation (\ref{l1}) is an
equality. Assume now that the right-hand side of equation (\ref{l1})
does not vanish. Similarly, we prove that in this case equation
(\ref{l1}) is an equality.

\medskip

\noindent\textbf{Lemma 3.} \textit{Let $q=4m+1$ where $m\not\in \{0,
-1\}$. Let $|2m+1|=p_1^{l_1}\times\cdots\times p_k^{l_k}$ --- be a
decomposition of $|2m+1|$ into prime factors. Let
$X=\Delta_{p_1}\times\cdots\times \Delta_{p_k}$. Then there exist
independent identically distributed random variables $\xi_1,\xi_2$
with values in $X$ and distribution $\mu\not\in I(X)$ such that the
conditional distribution of the linear form $L_2=\xi_1+q\xi_2$ given
$L_1=\xi_1+\xi_2$ is symmetric.}

\medskip

\textbf{Proof.} Since $X=\Delta_{p_1}\times\cdots\times
\Delta_{p_k}$, we have $Y\approx
\mathbb{Z}({p_1}^{\infty})\times\cdots\times
\mathbb{Z}({p_k}^{\infty})$. To avoid introducing new notation we
will assume that $Y=\mathbb{Z}({p_1}^{\infty})\times...\times
\mathbb{Z}({p_k}^{\infty})$.

Let $g_0(y)$ be an arbitrary characteristic functions on
$Y_{(2m+1)}$. Set

$$g(y)=
\left\{%
\begin{array}{ll}
    g_0(y), & \hbox{$y\in Y_{(2m+1)}$;} \\
    0, & \hbox{$y\not\in Y_{(2m+1)}$.} \\
\end{array}%
\right.
$$
The function $g(y)$ is a positive definite function on $Y$ (\cite[\S
32]{HeRo2}]). By the Bochner theorem there exists a distribution
$\mu \in {M}^1(X)$ such that $\widehat\mu(y) = g(y)$, $j =1, 2$. It
is clear that $g_0(y)$ can be chosen in such a way that $\mu\not\in
I(X)$.

Let $\xi_1,\xi_2$ be independent identically distributed random
variables with values in $X$ and   distribution $\mu$. We check that
the conditional distribution of the linear form $L_2=\xi_1+q\xi_2$
given $L_1=\xi_1+\xi_2$ is symmetric. By Lemma 1 it suffices to show
that the characteristic function $\hat\mu(y)$ satisfies equation
(\ref{l0}) which takes the form (\ref{l1}).

Let $u, v \in Y_{(2m+1)}$. Then $qv=(q+1)v-v=-v$ and equation
(\ref{l1}) is an equality.

If either $u\in Y_{(2m+1)}, v\not\in Y_{(2m+1)}$ or $u\not\in
Y_{(2m+1)}, v\in Y_{(2m+1)}$, then $u\pm v\not\in Y_{(2m+1)}$. Hence
$\hat\mu(u+v)=\hat\mu(u-v)=0$ and equation (\ref{l1}) is an
equality.

Let $u, v\not\in Y_{(2m+1)}$. Suppose that the left-hand side of
equation (\ref{l1}) does not vanish. Then $u+v\in Y_{(2m+1)}$,
$u+qv\in Y_{(2m+1)}$. Hence $(q-1)v \in Y_{(2m+1)}$. Since $q-1=4m$
and $f_{4m}\in\mathrm{Aut}(Y)$, we have $v\in Y_{(2m+1)}$, contrary
to the choice of $v$. Hence the left-hand side of equation
(\ref{l1}) is equal to zero. Similarly, we prove that the right-hand
side of equation (\ref{l1}) is equal to zero.

\medskip

\noindent\textbf{Lemma 4.} \textit{Let $X=\Sigma_{\mathbf{a}}$. If
$f_n\in \mathrm{Aut}(X)$, where $n=p_1^{l_1}\times\cdots\times
p_k^{l_k}$ is a decomposition of $n$ into prime factors, then the
group $X$ contains a subgroup topologically isomorphic to
$\Delta_{p_1}\times\cdots\times\Delta_{p_k}$.}

\medskip

\textbf{Proof.} Since $X=\Sigma_{\mathbf{a}}$, the character group
$Y=H_{\mathbf{a}}$ is a subgroup of $\mathbb{Q}$. As is well known
that
$$\mathbb{Q}/\mathbb{Z}\approx\mathop{\mbox{\rm\bf P}^*}\limits_{p \in
{\mathcal P}}\mathbb{Z}(p^{\infty}),$$ where ${\mathcal P}$ is the
set of prime numbers (\cite[\S 8]{Fu1}). Since $Y\subset
\mathbb{Q}$, we have $Y/\mathbb{Z} \subset \mathbb{Q}/\mathbb{Z}$.
The condition $f_n\in \mathrm{Aut}(X)$ implies that all $f_{p_j}\in
\mathrm{Aut}(X)$. Hence $f_{p_j}\in \mathrm{Aut}(Y)$. It is obvious
that if $p$ is a prime number and $f_{p}\in \mathrm{Aut}(Y)$ then
$F_p\subset Y/\mathbb{Z}$, where $F_p\approx
\mathbb{Z}(p^{\infty})$. Hence $L\subset Y/\mathbb{Z}$, where
$L\approx\mathbb{Z}(p_1^{\infty})\times\cdots\times
\mathbb{Z}(p_k^{\infty})$. It is clear that $Y/\mathbb{Z}=L\times
M$, where $M$ is a group. Since $(Y/\mathbb{Z})^*\approx
A(X,\mathbb{Z})\subset X$ and $(Y/\mathbb{Z})^*\approx L^*\times
M^*$, the group $X$ contains a subgroup topologically isomorphic to
$L^*\times M^*$. The statement of Lemma 4 follows from the form of
$L$.

\medskip

\noindent\textbf{Lemma 5.} \textit{Let $X$ be a locally compact
second countable Abelian group. Let $\xi_1, \ \xi_2$ be independent
random variables with values in $X$ and distributions $\mu_1,
\mu_2$. Consider the linear forms $L_1=\alpha_1\xi_1 +
\alpha_2\xi_2$ and $L_2=\beta_1\xi_1 + \beta_2\xi_2$, where
$\alpha_j, \beta_j$ are continuous homomorphisms of the group $X$.
The linear forms $L_1$ and $L_2$ are independent if and only if the
characteristic functions of the distributions $\mu_j$ satisfy the
equation
\begin{equation}\label{l5-1}
    \hat\mu_1(\tilde{\alpha_1} u+\tilde{\beta_1}v)
    \hat\mu_2(\tilde{\alpha_2} u+\tilde{\beta_2}v)=
    \hat\mu_1(\tilde{\alpha_1} u)
    \hat\mu_1(\tilde{\beta_1}v)
    \hat\mu_2(\tilde{\alpha_2}u)
    \hat\mu_2(\tilde{\beta_2}v),\quad u,v \in Y.
\end{equation}}

\medskip

Lemma 5 was proved in \cite[\S 10.1]{Fe-book2} in the case where
$\alpha_j, \beta_j \in {\rm Aut}(X)$. This proof is valid for
arbitrary continuous homomorphisms $\alpha_j, \beta_j$ of the group
$X$.

\medskip

\noindent\textbf{Lemma 6.} \textit{Let $X$ be a locally compact
second countable Abelian group, $\delta_1, \delta_2$ be continuous
homomorphisms of the group $X$. Let $\xi_1, \ \xi_2$ be independent
random variables with values in $X$ and distributions $\mu_1,
\mu_2$. If the conditional distribution of the linear form
$L_2=\delta_1\xi_1+\delta_2\xi_2$ given $L_1=\xi_1+\xi_2$ is
symmetric then the linear forms
$L'_1=(\delta_1+\delta_2)\xi_1+2\delta_2\xi_2$ and
$L'_2=2\delta_1\xi_1+(\delta_1+\delta_2)\xi_2$ are independent.}

\medskip

\textbf{Proof.} By Lemma 1 the symmetry of the conditional
distribution of the linear form  $L_2$ given $L_1$ implies that the
characteristic functions $\hat\mu_j(y)$ satisfy equation

\begin{equation}\label{l4.1}
    \hat\mu_1(u+\varepsilon_1 v)\hat\mu_2(u+\varepsilon_2 v)
    =\hat\mu_1(u-\varepsilon_1 v)\hat\mu_2(u-\varepsilon_2 v),\quad u,v\in Y,
\end{equation}
where $\varepsilon_j=\tilde{\delta_j}$.

Putting $u=\varepsilon_2 y, v=-y$ and then $u=-\varepsilon_1 y, v=y$
into equation (\ref{l4.1}) we obtain

\begin{equation}\label{l4.2}
    \hat\mu_1((\varepsilon_2 - \varepsilon_1)y)=
    \hat\mu_1((\varepsilon_1+\varepsilon_2)y)\hat\mu_2(2\varepsilon_2y),\quad y\in Y,
\end{equation}

\begin{equation}\label{l4.3}
    \hat\mu_2((\varepsilon_2-\varepsilon_1)y)=
    \hat\mu_1(-2\varepsilon_1y)\hat\mu_2(-(\varepsilon_1+\varepsilon_2)y),\quad y\in Y.
\end{equation}

Let $t, s\in Y$. Putting $u=\varepsilon_1s+\varepsilon_2t$, $v=s+t$
into (\ref{l4.1}) we obtain

\begin{equation}\label{l4.5}
    \hat\mu_1((\varepsilon_1+\varepsilon_2)t+2\varepsilon_1s)
    \hat\mu_2(2\varepsilon_2t+(\varepsilon_1+\varepsilon_2)s)
    =\hat\mu_1((\varepsilon_2-\varepsilon_1)t)
    \hat\mu_2(-(\varepsilon_2-\varepsilon_1)s),\quad s,t\in Y.
\end{equation}

Taking into account (\ref{l4.2}) and (\ref{l4.3}) equation
(\ref{l4.5}) can be written in the form

\begin{equation}\label{l4.6}
    \hat\mu_1((\varepsilon_1+\varepsilon_2)t+2\varepsilon_1s)
    \hat\mu_2(2\varepsilon_2t+(\varepsilon_1+\varepsilon_2)s)=$$ $$
    \hat\mu_1((\varepsilon_1+\varepsilon_2)t)\hat\mu_2(2\varepsilon_2t)
    \hat\mu_1(2\varepsilon_1s)\hat\mu_2((\varepsilon_1+\varepsilon_2)s),\quad s,t\in
    Y.
\end{equation}

Lemma 5 and equation (\ref{l4.6}) imply that the linear forms
$L'_1=(\delta_1+\delta_2)\xi_1+2\delta_2\xi_2$ and
$L'_2=2\delta_1\xi_1+(\delta_1+\delta_2)\xi_2$ are independent.

\medskip

\textbf{Remark 1.} Lemma 6 implies that the Heyde theorem on the
group $\mathbb{R}$ for $n=2$ can be obtained from the
Skitovich-Darmois theorem.

\medskip

\textbf{Proof of Theorem 1.} By Lemma 1 the symmetry of the
conditional distribution of the linear form $L_2$ given $L_1$
implies that the characteristic functions of distributions $\mu_j$
satisfy equation (\ref{l0}) which takes the form

\begin{equation}\label{1}
    \hat\mu_1(u+pv)\hat\mu_2(u+qv)=
    \hat\mu_1(u-pv)\hat\mu_2(u-qv),\quad u,v\in Y.
\end{equation}
We will study the solutions of this equation.

Consider first the case where $pq=-3$. Obviously, without loss of
generality we can assume that $p=1$ and $q=-3$ that is $L_1 = \xi_1
+ \xi_2$ and $L_2 = \xi_1 -3\xi_2$. Lemma 6 implies that the linear
forms $L'_1=-2\xi_1-6\xi_2$ and $L'_2=2\xi_1-2\xi_2$ are
independent. Making the substitution $\zeta_1=2\xi_1$ and
$\zeta_2=-2\xi_2$, we obtain that the linear forms
$L''_1=-\zeta_1+3\zeta_2$ and $L''_2=\zeta_1+\zeta_2$ are also
independent. As has been proved in \cite{Fe-solenoid-2009} the
independence of the linear forms $L''_1$ and $L''_2$ implies that at
least the distribution of one random variable $\zeta_j$ belongs to
$\Gamma(X)*I(X)$. Returning to the random variables $\xi_j$, we
obtain the statement 1 of Theorem 1.

Consider now the case where $pq\neq -3$. Two cases are possible:
$pq$ is a composite number and $pq$ is a prime number.

We prove that in these cases there exist independent random
variables $\xi_1$ and $\xi_2$ with values in $X$ and distributions
$\mu_1, \mu_2 \notin \Gamma (X)*I(X)$ such that the conditional
distribution of the linear form $L_2$ given $L_1$ is symmetric.

\medskip

\textbf{Case} $\mathbf{1.}$ $pq$ is a composite number. In this case
we follow the scheme of the proof of the analogous case in Theorem 1
of the article \cite{Fe-solenoid-2009}.

\medskip

Put $s=p-q$, and decompose $|s|$ into prime factors $|s|=s_1^{k_1}
\cdots s_l^{k_l}$. Denote by $H$ the subgroup of $Y$ of the form

$$ H = \left\{{m \over s_{j_1}^{n_1} \cdots s_{j_r}^{n_r}} : \quad
m, n_j \in {\mathbb{Z}} \right\}. $$ If $|s|=1$ we suppose that
$H=\mathbb{Z}$. Set $G=H^*$.

\medskip

$\mathbf{1a.}$ $|p|>1, |q|>1$.

\medskip

Since $p$ and $s$ are relatively prime, and so $q$ and $s$, we have
$H^{(p)} \ne H$ and $H^{(q)} \ne H$. Assume that $\lambda_j \in
{M}^1(G)$ and $\sigma(\lambda_1) \subset A(G, H^{(p)}), \
\sigma(\lambda_2) \subset A(G, H^{(q)})$. It follows from this that
$\widehat\lambda_1(y)=1$, $y\in H^{(p)}$, and
$\widehat\lambda_2(y)=1$, $y\in H^{(q)}$. Therefore

\begin{equation}
\label{t1.1} \widehat\lambda_1(u + pv) = ^{}\widehat\lambda_1(u), \
\widehat\lambda_2(u + qv) = \widehat\lambda_2(u), \quad u, v \in H.
\end{equation}

\noindent Consider the functions $g_j(y)$ on the group $Y$ of the
form

\begin{equation}
\label{t1.2} g_j(y) =
\left\{%
\begin{array}{ll}
    \widehat\lambda_j(y), & \hbox{$y \in H$;} \\
    0, & \hbox{$y \notin H$.} \\
\end{array}%
\right.
\end{equation}

\noindent The functions $g_j(y)$ are positive definite functions on
$Y$ (\cite[\S32]{HeRo2}). By the Bochner theorem there exist
distributions $\mu_j \in {M}^1(X)$ such that $\widehat\mu_j(y) =
g_j(y)$, $j =1, 2$. We will show that the characteristic functions
$\hat\mu_j(y)$ satisfy equation (\ref{1}).

We conclude from (\ref{t1.1}) and (\ref{t1.2}) that if $u, v \in H$,
then equation (\ref{1}) is an equality.

Let either $u\in H, v\not\in H$ or $u\not\in H, v\in H$. Since the
numbers $p$ and $s$ are relatively prime, we have either $pv\not\in
H$ or $pv\in H$ respectively. So, $u\pm pv\not\in H$ and hence
$\hat\mu_1(u+pv)=\hat\mu_1(u-pv)=0$, and equation (\ref{1}) is an
equality.

Let $u, v\not\in H$. Suppose that the left-hand side of equation
(\ref{1}) does not vanish. Then $u+pv\in H$ and $u+qv\in H$. Hence
$sv \in H$. Therefore $v\in H$, contrary to the choice of $v$. Hence
the left-hand side of (\ref{1}) is equal to zero. Similarly, we
prove that the right-hand side of (\ref{1}) is equal to zero.

So, the characteristic functions $\widehat\mu_j(y)$ satisfy equation
(\ref{1}). If $\xi_1$ and $\xi_2$ are independent random variables
with values in $X$ and distributions $\mu_j$, then by Lemma 1 the
conditional distribution of the linear form $L_2$ given $L_1$ is
symmetric. It is clear that $\lambda_j$ can be chosen in such a way
that $\mu_1, \mu_2 \notin \Gamma (X)*I(X)$. The desired statement in
case \textbf{1a} is proven.

\medskip

$\mathbf{1b.}$ Either $|p|=1, |q|>1$ or $|p|>1, |q|=1$.

\medskip

Assume for definiteness that $|p|=1$. Without loss of generality, we
suppose $p=1$. Let $q = q_1 q_2$ be a decomposition of $q$, where
$|q_j|>1$, $j=1, 2$. It is obvious that if $f_q \in {\rm Aut}(X)$,
then $f_{q_1}, f_{q_{2}} \in {\rm Aut}(X)$. Note that the
conditional distribution of $L_2 = \xi_1 + q\xi_2$ given $L_1 =
\xi_1 + \xi_2$ is symmetric if and only if the conditional
distribution of $L_2 = {1 \over q_1}\xi_1 + q_2\xi_2$ given $L_1 =
\xi_1 + \xi_2$ is symmetric. Making the substitution $\zeta_1= {1
\over q_1} \xi_1$, we reduce the problem to the case when $L_1 =q_1
\xi_1 + \xi_2$, $L_2 = \xi_1 + q_2 \xi_2$. Equation (\ref{1}) in
this case takes the form

\begin{equation}
    \label{2} \hat\mu_1(q_1u + v)\hat\mu_2(u + q_2v)=
    \hat\mu_1(q_1u -v)\hat\mu_2(u - q_2v), \quad u, v \in Y.
\end{equation}

Assume that $\lambda_j \in {M}^1(G)$ and $\sigma(\lambda_j) \subset
A(G, H^{(q_j)})$, $j = 1, 2$. It is obvious that
$\widehat\lambda_j(y)=1$, $y\in H^{(q_j)}$. Hence
\begin{equation}
\label{t1.3}\widehat\lambda_1(q_1u + v) = \widehat\lambda_1(v),
\quad  \widehat\lambda_2(u + q_2v) = \widehat\lambda_2(u) \quad u, v
\in H.
\end{equation}
In the same manner as in case \textbf{1a} we define the functions
$g_j(y)$ by formulas (\ref{t1.2}) and distributions $\mu_j \in
{M}^1(X)$. We will show that the characteristic functions
$\hat\mu_j(y)$ satisfy equation (\ref{2}).

We conclude from (\ref{t1.3}) and (\ref{t1.2}) that if $u, v \in H$,
than equation (\ref{2}) is an equality.

Let either $u\in H, v\not\in H$ or $u\not\in H, v\in H$. Since the
numbers $q_1$ and $s$ are relatively prime, we have either $pu\in H$
or $pu\not\in H$ respectively. So, $q_1u\pm v\not\in H$ and hence
$\hat\mu_1(q_1u+v)=\hat\mu_1(q_1u-v)=0$, and equation (\ref{2}) is
an equality.

Let $u, v\not\in H$. Suppose that the left-hand side of equation
(\ref{2}) does not vanish. Then $q_1u+v\in H$ and $u+q_2v\in H$.
Hence $su \in H$. Therefore $u\in H$, contrary to the choice of $u$.
Hence the left-hand side of equation (\ref{2}) is equal to zero.
Reasoning similarly we show that the right-hand side of (\ref{2}) is
equal to zero.

So, the characteristic functions $\widehat\mu_j(y)$ satisfy equation
(\ref{2}). If $\xi_1$ and $\xi_2$ are independent random variables
with values in $X$ and distributions $\mu_j$, then by Lemma 1 the
conditional distribution of the linear form $L_2$ given $L_1$ is
symmetric. It is clear that $\lambda_j$ can be chosen in such a way
that $\mu_1, \mu_2 \notin \Gamma (X)*I(X)$. The desired statement in
case \textbf{1b} is proven.

\medskip

\textbf{Case} $\mathbf{2.}$ $pq$ is a prime number, i.e. either
$|p|=1, |q|>1$ or $|p|>1, |q|=1$.

\medskip

Assume for definiteness that $p=1$ and $q$ is a prime number, i.e.
$L_1 = \xi_1 + \xi_2$, $L_2 = \xi_1 + q\xi_2$. Equation (\ref{1})
takes the form

\begin{equation}\label{3}
    \hat\mu_1(u+v)\hat\mu_2(u+qv)=
    \hat\mu_1(u-v)\hat\mu_2(u-qv),\quad u,v \in Y.
\end{equation}

\medskip

$\mathbf{2a.}$ $|q|=2$.

\medskip

Since $f_2\in \mathrm{Aut}(X)$, Lemma 4 implies that the group $X$
contains a subgroup topologically isomorphic to $\Delta_2$. Then the
statement 2 of Theorem 1 follows from Lemma 2.

\medskip

Let $q$ be an odd number. There exist two possibilities: 1)
$q=4m+3$; 2) $q=4m+1$.

\medskip

$\mathbf{2b.}$ $q=4m+3$.

\medskip

Note that since $q$ is an odd number and $f_{q+1}\in \mathrm{Aut}
(X)$, the homomorphism $f_2\in \mathrm{Aut}(X)$. Hence Lemma 4
implies that the group $X$ contains a subgroup topologically
isomorphic to $\Delta_2$. Then the statement 2 of Theorem 1 follows
from Lemma 2.

\medskip

$\mathbf{2c.}$ $q=4m+1$ ($m\neq -1$).

\medskip

Since $f_{q+1}\in \mathrm{Aut}(X)$ and $q+1=2(2m+1)$, the
homomorphism $f_{2m+1}\in \mathrm{Aut}(X)$. Let
$|2m+1|=p_1^{l_1}\times...\times p_k^{l_k}$ be a decomposition of
$|2m+1|$ into prime factors. Lemma 4 implies that the group $X$
contains a subgroup topologically isomorphic to
$\Delta_{p_1}\times...\times\Delta_{p_k}$. Then the statement 2 of
Theorem 1 follows from Lemma 3.

\medskip

\textbf{Remark 2.} Statement 1 of Theorem 1 may not be strengthened.
Namely, the following statement is valid. If $pq=-3$ then there
exist independent random variables $\xi_1,\xi_2$ with values in $X$
and distributions $\mu_1,\mu_2$ such that the conditional
distribution of the linear form $L_2=p\xi_1+q\xi_2$ given
$L_1=\xi_1+\xi_2$ is symmetric and one of the distributions
$\mu_j\not\in\Gamma(X)*I(X)$.

It is suffices to consider the case when $p=1$, $q=-3$. We shall
construct solutions of equation (\ref{l0}) which takes the form

\begin{equation}\label{r1.2}
    \hat\mu_1(u+v)\hat\mu_2(u-3v)=\hat\mu_1(u-v)\hat\mu_2(u+3v),\quad u,v\in Y.
\end{equation}

Let $\gamma_1$ and $\gamma_2$ be Gaussian distributions on $X$ with
characteristic functions $\hat\gamma_1(y)=e^{-3y^2}$ and
$\hat\gamma_2(y)=e^{-y^2}$. It is easy to verify that these
functions satisfy equation (\ref{r1.2}).

Since $f_{p+q}\in \mathrm{Aut}(X)$, we have $f_{2}\in
\mathrm{Aut}(X)$. Hence $f_{2}\in \mathrm{Aut}(Y)$. Therefore the
group $Y$ contains a subgroup of dyadic rational numbers. Denote by
$H$ this subgroup. Since $f_{q}\in \mathrm{Aut}(X)$, we have
$f_{3}\in \mathrm{Aut}(X)$. Hence $f_{3}\in \mathrm{Aut}(Y)$.
Therefore $1/3 \in Y$. Denote by $L$ a subgroup in $Y$ which is
generated by the subgroup $H$ and the element $1/3$. Observe that
$L=\{H, 1/3+H, 2/3+H\}$. Let $G=A(X,H)$, $K=A(X,L)$. Let
$\omega_1=(1/2)[ m_G+  m_K]$ and $\omega_2=m_G$. It follows from
(\ref{def1}) that

\begin{equation}\label{r1.1}
    \hat\omega_1(y)=
        \left\{%
\begin{array}{ll}
    1, & \hbox{$y \in H$;} \\
    1/2, & \hbox{$y \in L\setminus H$;} \\
    0, & \hbox{$y \notin L$.} \\
\end{array}%
\right.
    \quad
    \hat\omega_2(y)=
\left\{%
\begin{array}{ll}
    1, & \hbox{$y \in H$;} \\
    0, & \hbox{$y \notin H$.} \\
\end{array}%
\right.
\end{equation}
We verify that these functions satisfy equation (\ref{r1.1}).

It is clear that if $u\in H$, $v\in L$ then equation (\ref{r1.1}) is
an equality.

If $u\in L\setminus H$, $v\in L$ then $u\pm 3v\not\in H$. Hence
$\hat\omega_2(u-3v)=\hat\omega_2(u+3v)=0$ and equation (\ref{r1.1})
is an equality.

If either $u\in L$, $v\not\in L$ or $u\not\in L$, $v\in L$ then
$u\pm v\not\in L$. Hence $\hat\omega_1(u+v)=\hat\omega_2(u-v)=0$ and
equation (\ref{r1.1}) is an equality.

Let $u, v\not\in L$. Suppose that the left-hand side of equation
(\ref{r1.1}) does not vanish. Then $u+v\in L$ and $u-3v\in H$. Hence
$4v \in L$. We obtain that $v\in L$, contrary to the choice of $v$.
Hence the left-hand side of equation (\ref{r1.1}) is equal to zero.
Reasoning similarly we show that the right-hand side of equation
(\ref{r1.1}) is equal to zero.

Put $\mu_j=\gamma_j*\omega_j$, $j=1,2$. It is obvious that the
functions $\hat\mu_j(y)$ satisfy equation (\ref{r1.1}). If $\xi_1$
and $\xi_2$ are independent random variables with values in $X$ and
distributions $\mu_j$ then Lemma 1 implies that the conditional
distribution of the linear form $L_2=\xi_1-3\xi_2$ given
$L_1=\xi_1+\xi_2$ is symmetric. By the construction
$\mu_1\not\in\Gamma(X)*I(X)$ and $\mu_2\in\Gamma(X)*I(X)$.

\medskip

\textbf{Remark 3.} Note that in Theorem 1 we suppose that there
exist for some relatively prime $p$ and $q$ automorphisms $f_p$ and
$f_q$ such that $f_{p\pm q}\in \mathrm{Aut}(X)$ on the group
$X=\Sigma_{\mathbf{a}}$. It is easy to prove that groups
$X=\Sigma_{\mathbf{a}}$ have this property if and only if $f_2, f_3
\in \mathrm{Aut}(X)$.

\medskip

\textbf{Remark 4.} We note that if in Theorem 1 distributions
$\mu_1, \mu_2$ have non-vanishing characteristic functions,  then
$\mu_1, \mu_2 \in \Gamma(X)$. Indeed, it follows from conditions on
coefficients of the linear forms that one of numbers $p, q, p\pm q$
is even. So, $f_2\in \mathrm{Aut}(X)$. Hence the group
$X=\Sigma_{\mathbf{a}}$ does not contain elements of order two. The
following theorem (see \cite{Fe-He-2}) implies the desired
statement: Let $X$ be a locally compact second countable Abelian
group containing no elements of order two. Let $\xi_1,\xi_2$ be
independent random variables with values in $X$ and distributions
$\mu_1,\mu_2$ with nonvanisging characteristic functions. Consider
the linear forms $L_1=\xi_1 + \xi_2$ and $L_2=\delta_1\xi_1 +
\delta_2\xi_2$, where $\delta_j, \delta_1\pm \delta_2 \in {\rm
Aut}(X)$. If the conditional distribution of the linear form $L_2$
given $L_1$ is symmetric then $\mu_1, \mu_2\in \Gamma(X)$.

\section{The Heyde theorem (the special case)}

We prove in this section that Theorem 1 can be essentially
strengthened if we assume in addition  that the support at least one
of   distributions $\mu_j$ is not contained in a coset of a proper
closed subgroup of the group $X$.

Let $\mu\in M^1(X)$. It is easy to see that  $\mu$  has the
property: $\sigma(\mu)$ is not contained in a coset of a proper
closed subgroup of the group $X$ if and only if

\begin{equation}\label{intro2.1}
    \{ y\in Y:\ |\hat\mu(y)|=1\}=\{0\}.
\end{equation}

\medskip

\textbf{Theorem 2.} \textit{Let $X=\Sigma_{\mathbf{a}}$. Assume that
$f_p, f_q, f_{p\pm q} \in {\rm Aut}(X)$, $p$ and $q$ are relatively
prime. The following statements hold:}
\begin{enumerate}
    \item \textit{Let $pq>0$. Let $\xi_1,\xi_2$
          be independent random variables with values in $X$
          and distributions $\mu_1,\mu_2$ such that at least
          one   support  $\sigma(\mu_j)$ is not contained in a
          coset of a proper closed subgroup of the group $X$.
          If the conditional distribution of the linear form
          $L_2=p\xi_1+q\xi_2$ given $L_1=\xi_1+\xi_2$
          is symmetric then $\mu_1=\mu_2=m_X$.}
    \item \textit{Let $pq= -3$. Let $\xi_1,\xi_2$
          be independent random variables with values in $X$
          and distributions $\mu_1,\mu_2$ such that at least
          one    support  $\sigma(\mu_j)$ is not contained in a
          coset of a proper closed subgroup of the group $X$.
          If the conditional distribution of the linear form
          $L_2=p\xi_1+q\xi_2$ given $L_1=\xi_1+\xi_2$
          is symmetric then at least one distribution
          $\mu_j\in \Gamma(X)*I(X)$.}
    \item \textit{Let $pq<0$ and $pq\neq -3$. Then there exist
          independent random variables $\xi_1,\xi_2$
          with values in $X$ and distributions
          $\mu_1,\mu_2$ such that the conditional
          distribution of the linear form $L_2=p\xi_1+q\xi_2$
          given $L_1=\xi_1+\xi_2$ is symmetric, the
          distributions $\mu_j\not\in \Gamma(X)*I(X)$,
          and each of the supports $\sigma(\mu_j)$
          is not contained in a coset of a proper closed
          subgroup of the group $X$.}
\end{enumerate}

\medskip

To prove Theorem 2 we need some lemmas. The lemmas 7 and 8 given
below were proved in \cite{Fe-solenoid-2009} in the case $a=1$. In
the case $a\neq 1$ proofs of lemmas 7 and 8 follow the schemes of
proofs of the corresponding lemmas in \cite{Fe-solenoid-2009}.

\medskip

\textbf{Lemma 7.} \textit{Let $Y$ be an arbitrary Abelian group, let
$a, b \in {\mathbb{Z}}$, $ab\ne 0$, and let $g_1(y)$ and $g_2(y)$ be
functions on $Y$ satisfying the equation
\begin{equation}\label{l5.1}
    g_1(u+av)g_2(u+bv)=g_1(u)g_1(av)g_2(u)g_2(bv),
    \quad u, v \in Y,
\end{equation}
and the conditions
\begin{equation}\label{l5.1.1}
 g_1(-y)=\overline{g_1(y)}, \quad g_2(-y)=\overline{g_2(y)}, \quad y \in Y, \quad g_1(0)=g_2(0)=1.
\end{equation}
Set $c=a-b$. If for certain $y_0 \in Y^{(c)}$ the inequality
$g_1(y_0)g_2(y_0) \ne 0$ holds thenthere exists a subgroup
$M=\{kabz_0\}_{k \in {\mathbb{Z}}}$ ($y_0=cz_0$, $z_0\in Y$) such
that $g_1(y)g_2(y) \ne 0$ for $y \in M$.}

\medskip

\textbf{Proof.} Putting $u = -by, \ v = y$ and then $u = ay, \ v =
-y$ in (\ref{l5.1}) we get

\begin{equation} \label{l5.2}
    g_1(cy)=g_1(-by)g_1(ay)g_2(-by)g_2(by), \quad y \in Y,
\end{equation}

\begin{equation} \label{l5.3}
    g_2(cy)=g_1(ay)g_1(-ay)g_2(ay)g_2(-by), \quad y \in Y.
\end{equation}

By the condition of the lemma $y_0=cz_0$, $z_0\in Y$. Substituting
$y=z_0$ into (\ref{l5.2}) and (\ref{l5.3}) we conclude that

\begin{equation} \label{l5.4}
    g_1(az_0) \ne 0, g_1(bz_0) \ne 0, \quad
    g_2(az_0) \ne 0, g_2(bz_0) \ne 0.
\end{equation}

Putting $u=az_0, \ v=kz_0$ and then $u=bz_0, \ v=kz_0, \ k \in
{\mathbb{Z}}$, in equation (\ref{l5.1}) we obtain

\begin{equation} \label{l5.5}
    g_1((k+1)az_0)g_2((bk+a)z_0)=g_1(az_0)g_1(akz_0)g_2(az_0)g_2(bkz_0),
\end{equation}

\begin{equation} \label{l5.6}
    g_1((ak+b)z_0)g_2((k+1)bz_0)=g_1(bz_0)g_1(akz_0)g_2(bz_0)g_2(bkz_0).
\end{equation}

Taking into account (\ref{l5.4}), it follows by induction from
(\ref{l5.5}) and (\ref{l5.6}) that $g_1(kaz_0) \ne 0, \ g_2(kbz_0)
\ne 0, \ k \in {\mathbb{Z}}$. The subgroup $M=\{kabz_0\}_{k \in
{\mathbb{Z}}}$ is the required one.

\medskip

\textbf{Lemma 8.} \textit{Let $M$ be an arbitrary subgroup in
$\mathbb{Q}$, $g_1(y)$ and $g_2(y)$ be functions on $M$ satisfying
equation $(\ref{l5.1})$, conditions $(\ref{l5.1.1})$, and the
conditions
\begin{equation} \label{l7.1}
    0 < g_1(y) \le 1, \quad 0 < g_2(y) \le 1.
\end{equation}
Put $c=b-a$. Then on the subgroup $M^{(cab)}$ the following
representation holds:
\begin{equation} \label{l7.1.1}
    g_2(y)=\exp\{-\lambda_1y^2\}, \quad g_2(y)=\exp\{- \lambda_2y^2\},
\end{equation}
where $\lambda_j \ge 0$.}

\medskip

\textbf{Proof.} Set $\varphi_1(y)= - \ln g_1(y), \ \varphi_2(y)= -
\ln g_2(y)$. It follows from (\ref{l5.1}) that

\begin{equation}\label{l7.2}
    \varphi_1(u+av) + \varphi_2(u+bv) = A(u) + B(v), \quad u, v \in M,
\end{equation}
where $A(u) = \varphi_1(u) + \varphi_2(u)$, $B(v) = \varphi_1(av) +
\varphi_2(bv)$.

We use the finite difference method to solve equation (\ref{l7.2}).

Let $h_1$ be an arbitrary element of $M$. Substitute $u+bh_1$ for
$u$ and $v-h_1$ for $v$ in equation (\ref{l7.2}) and subtract
equation (\ref{l7.2}) from the resulting equation. We get

\begin{equation}\label{l7.3}
    \Delta_{ch_1} \varphi_1(u+av) = \Delta_{bh_1}A(u) +
    \Delta_{-h_1}B(v).
\end{equation}
Putting $v=0$ in (\ref{l7.3}) and subtracting the resulting equation
from (\ref{l7.3}) we obtain

\begin{equation}\label{l7.4}
    \Delta_{av}\Delta_{ch_1} \varphi_1(u) =\Delta_{-h_1}B(v)-\Delta_{-h_1}B(0).
\end{equation}
Substitute $u+h_1$ for $u$ in equation (\ref{l7.4}) and subtract
equation (\ref{l7.4}) from the resulting equation. We get

\begin{equation}\label{l7.5}
    \Delta_{h_1}\Delta_{av}\Delta_{ch_1} \varphi_1(u) =0.
\end{equation}
We conclude from (\ref{l7.5}) that the function $\varphi_1(y)$
satisfies the equation

\begin{equation}\label{l7.6}
    \Delta_{h}^3 \varphi_1(y) =0, \quad y\in M, h\in M^{(ca)}.
\end{equation}

Reasoning similarly we get

\begin{equation}\label{l7.7}
    \Delta_{h}^3 \varphi_2(y) =0, \quad y\in M, h\in M^{(cb)}.
\end{equation}

It follows from (\ref{l7.6}) and (\ref{l7.7}) that the functions
$\varphi_j(y)$ are polynomials of the degree 2 on the subgroup
$M^{(cab)}$. Taking into account (\ref{l5.1.1}) and (\ref{l7.1}) on
the subgroup $M^{(cab)}$ we get $\varphi_j(y)=\lambda_j y^2$ where
$\lambda_j \ge 0$.

\medskip

\textbf{Proof of Theorem 2.} Let $pq>0$. Lemma 6 implies that the
linear forms $L'_1=(p+q)\xi_1+2q\xi_2$ and $L'_2=2p\xi_1+(p+q)\xi_2$
are independent. Making the substitution $\xi'_1=(p+q)\xi_1$ and
$\xi'_2=2q\xi_2$, we obtain that the linear forms
$L'_1=\xi'_1+\xi'_2$ and $L'_2={2p\over p+q}\xi'_1+{p+q\over
2q}\xi'_2$ are also independent. We also note that if ${\delta \in
{\rm Aut}(X)}$ then the linear forms $L_1$ and $L_2$ are independent
if and only if the linear forms $L_1$ and $\delta L_2$ are
independent. Thus we may assume without loss of generality that
$L'_1=\xi'_1+\xi'_2$ and $L'_2=4pq\xi'_1+(p+q)^2\xi'_2$. Denote by
$\mu'_j$ the distributions of random variables $\xi'_j$. Since $f_2,
f_p, f_q, f_{p+ q}\in \mathrm{Aut}(X)$, if we prove that
$\mu'_j=m_X$ then Theorem 2 in case 1 will be proved.

By Lemma 5 the independence of $L'_1$ and $L'_2$ implies that the
characteristic functions of distributions $\mu'_j$ satisfy equation
(\ref{l5-1}) which takes the form

\begin{equation}\label{t2.1}
    \hat\mu'_1(u+4pqv)\hat\mu'_2(u+(p+q)^2v)=
    \hat\mu'_1(u)\hat\mu'_1(4pqv)
    \hat\mu'_2(u)\hat\mu'_2((p+q)^2v), \quad u,v\in Y.
\end{equation}

It is clear that the characteristic functions of distributions
$\bar\mu'_j$ also satisfy equation (\ref{t2.1}). Therefore the
characteristic functions of distributions $\nu_j=\mu'_j*\bar\mu'_j$
satisfy equation (\ref{t2.1}). Note that
$\widehat\nu_j(y)=|\widehat\mu'_j(y)|^2 \ge 0, \ j = 1, 2$. We also
note that since at least one support $\sigma(\mu_j)$ is not
contained in a coset of a proper closed subgroup of the group $X$,
we have that at least one support $\sigma(\nu_j)$ is not contained
in a coset of a proper closed subgroup of the group $X$. It follows
from (\ref{intro2.1}) that at least for one $j$ the equality

\begin{equation}\label{t2.4}
    \{y\in Y:\ \hat\nu_j(y)=1\}=\{0\}
\end{equation}
holds.

Putting $u=-(p+q)^2y, v=y$ and then $u=-4pqy, v=y$ into equation
(\ref{t2.1}) we obtain

\begin{equation}\label{t2.2}
    \hat\nu_1((p-q)^2y)=
    \hat\nu_1((p+q)^2y)\hat\nu_1(4pqy)
    \hat\nu_2^2((p+q)^2y), \quad y\in Y.
\end{equation}

\begin{equation}\label{t2.3}
    \hat\nu_2((p-q)^2y)=
    \hat\nu_1(4pqy)
    \hat\nu_2(4pqy)\hat\nu_2((p+q)^2y), \quad y\in Y.
\end{equation}

Assume first that $\hat\nu_1(y)\hat\nu_2(y) = 0$, $y \in Y$, $y\neq
0$. It follows from (\ref{t2.2}) that $\hat\nu_1((p-q)^2y) = 0$, $y
\in Y$, $y\neq 0$. Since $f_{p-q} \in {\rm Aut}(X)$, we conclude
that $\hat\nu_1(y) = 0$, $y \in Y$, $y\neq 0$. Hence $\nu_1 = m_X$,
so that $\mu'_1 = m_X$. Similarly, (\ref{t2.3}) implies that $\mu'_2
= m_X$.

Assume now that $\hat\nu_1(y_0)\hat\nu_2(y_0) \ne 0$ for some $y_0
\in Y, \ y_0 \ne 0$. Since $f_{p-q}\in \mathrm{Aut}(X)$, we have
that $Y^{((p-q)^2)} = Y$. We can apply Lemma 7 and obtain a subgroup
$M \subset Y$ such that $\hat\nu_1(y)\hat\nu_2(y) \ne 0$. By Lemma 8
the restrictions of the characteristic functions $\hat\nu_1(y)$ and
$\hat\nu_2(y)$ to $M^{(4pq(p+q)^2)(p-q)^2)}$ have form
(\ref{l7.1.1}). Substituting these representations into (\ref{t2.1})
we get

$$4pq\lambda_1+(p+q)^2\lambda_2=0.$$
Since $pq>0$, this equality implies that $\lambda_1=\lambda_2=0$.
Hence $\hat\nu_1(y)=\hat\nu_2(y)=1$ for $y\in
M^{(4pq(p+q)^2)(p-q)^2)}$, that contradicts (\ref{t2.4}).

Let $pq=-3$. The desired statement follows from statement 1 of
Theorem 1.

Let $pq<0$, $pq\neq -3$. Denote by $\omega_j$ the distributions
constructed in the proof of Theorem 1 in corresponded cases. We note
that the characteristic functions $\hat\omega_1(y)$ and
$\hat\omega_2(y)$ satisfy equation (\ref{1}). Denote by $\gamma_j$
Gaussian distributions on $X$ with the characteristic functions
$\hat\gamma_1(y)=e^{-\lambda y^2}$, $\hat\gamma_2(y)=e^{{p\over q}
\lambda y^2}$, where $\lambda>0$. It is easy to verify that the
functions $\hat\gamma_1(y)$ and $\hat\gamma_2(y)$ satisfy equation
(\ref{1}). Put $\mu_j=\omega_j*\gamma_j$. It is obvious that the
functions $\hat\mu_1(y)$ and $\hat\mu_2(y)$ satisfy equation
(\ref{1}). Since a support of a symmetric non degenerate Gausssian
distribution is a connected subgroup, we have that
$\sigma(\gamma_j)=X$. Hence $\sigma(\mu_j)=X$. By the construction
$\mu_j\not\in \Gamma(X)*I(X)$. Thus  $\mu_j$ are the desired
distributions.

\medskip

\textbf{Remark 5.} The example of distributions constructed in
Remark 2 shows that statement 2 of Theorem 2 may not be strengthened
to the statement that both distributions $\mu_j\in\Gamma(X)*I(X)$.

\medskip

\bibliographystyle{model1a-num-names}

\end{document}